\documentclass[11pt]{amsart}
\usepackage{latexsym,amssymb,amsmath}
\def\xxx#1#2{x^{#1}_{#2}}        
\def\t{\tau}
\def\s{\sigma}

\def\frp#1{\frac{\partial}{\partial{#1}}}

\def\BC{\mathbb C}

\def\BP{\mathbb P}
\def\pp#1{\mathbb P^{#1}}

\def\pp#1{{\mathbb P}^{#1}}
\def\tdim{\rm dim}
\def\hd{,...,}
\def\ww{\wedge}
\def\upperp{{}^\perp}

\def\ctimes{\times \cdots\times}
\def\tlker{\text{Lker}}\def\trker{\text{Rker}}

\def\11{\mathbf 1}

\def\FS{{\mathfrak S}}

\def\l{\lambda}
\def\a{\alpha}

\def\b{\beta}
\def\g{\gamma}
\def\s{\sigma}

\def\ot{{\mathord{\,\otimes }\,}}

\def\otc{{\mathord{\otimes\cdots\otimes}\;}}

\def\ra{{\mathord{\;\rightarrow\;}}}

\def\La#1{\Lambda^{#1}}
\def\tim{\text{Image}\,}

\def\tdim{\text{dim}\,}

\def\tmod{\text{ mod }}

\def\tlim{\text{lim}}

\def\wcdots{\ww \cdots\ww}

\newtheorem{theorem}{Theorem}[section]
\newtheorem{proposition}[theorem]{Proposition}
\newtheorem{lemma}[theorem]{Lemma}

\theoremstyle{definition}

\theoremstyle{remark}
\newtheorem{remark}[theorem]{Remark}

\begin{document}

\title{The border rank of the multiplication of $2\times 2$ matrices is seven}
\author{J.M. Landsberg}
\email{jml@math.tamu.edu}

\maketitle

\section{Introduction}
   One of the leading problems of 
 algebraic complexity theory is matrix multiplication. 
The na\"\i ve multiplication of two $n\times n$ matrices uses $n^3$ multiplications.
 In 1969 Strassen \cite{Strassen493}   presented
    an explicit algorithm for multiplying $2\times 2$ matrices
 using seven multiplications. 
 In the opposite direction, Hopcroft and Kerr \cite{hk252} and Winograd \cite{win553}
 proved independently that there is no algorithm for multiplying $2\times 2$ matrices using
 only $6$ multiplications.
 
 The precise number of multiplications needed to execute matrix multiplication 
(or any given bilinear map) is   called the {\it rank} of the bilinear map.
  A related problem is to determine the {\it border rank} of
 matrix multiplication (or any given bilinear map),
 first introduced in  \cite{MR592760,MR534068}.
Roughly speaking, some bilinear maps may be approximated
 with arbitrary precision by less complicated bilinear maps and the border rank
 of a bilinear map is the complexity of arbitrarily small \lq\lq good\rq\rq\  perturbations
 of   the map. These perturbed maps  
 can give rise to fast exact algorithms for matrix multiplication,
 see \cite{BCS}. The border rank made   appearances in the literature in
 the 1980's and early 90's, see, e.g.,
\cite{MR592760,MR534068,MR92m:11038,MR91i:68058,MR91f:68099,MR88g:68048,MR88d:68020,MR88g:15021, 
MR87j:15055,MR87g:65053,MR87f:15017,MR86c:68040,Strassen505,MR82e:68037,MR81i:68004,Strbook},
but to our  knowledge there has not been much progress on the question since
then.

\smallskip

More precisely, for any complex projective variety $X\subset \BC\BP ^N=\BP V$, and
point  $p\in \BP V$,
define the {\it $X$-rank of $p$} to be the smallest number $r$ such that $p$ is in the linear
span of $r$ points of $X$. Define
$\s_r(X)$, the {\it $r$-th secant variety of $X$},
 to be the Zariski closure of the set of points of $X$-rank $r$,
and define the {\it $X$-border rank of $p$} to be the smallest $r$ such that
$p\in \s_r(X)$. The terminology is motivated by the case $X=Seg(\pp {a-1}\times \pp {b-1})\subset
\BP (\BC^{a }\otimes \BC^{b })$, the {\it Segre variety} of rank one matrices. Then
the $X$-rank of a matrix is just its usual rank.

Let $A^*,B^*,C$ be vector spaces and let $f: A^*\times B^*\ra C$ be a bilinear map, i.e.,
an element of $A\ot B\ot C$. Let $X=Seg(\BP A\times \BP B\times \BP C)\subset \BP (A\ot B\ot C)$
denote the {\it Segre variety} of decomposable tensors in $A\ot B\ot C$.  The border rank of a bilinear map is its $X$-border rank.
While for the
Segre product of two projective spaces, border rank   coincides with rank, here
they can be quite different. 

\smallskip

In this paper we prove the theorem stated in the title. 
Let $MMult\in \BC^4\ot \BC^4\ot \BC^4$ denote the matrix multiplication operator
for two by two matrices.
Strassen \cite{Strassen505}   showed that $MMult\notin \s_5(\pp 3\times \pp 3\times \pp 3)$.
Our method of proof is to   decompose $\s_6(\pp 3\times \pp 3\times \pp 3)\backslash \s_5(\pp 3\times \pp 3\times \pp 3)$ into
 sixteen components, and then using case by case arguments show $MMult$
is not in any of the components. The decomposition rests upon a
differential-geometric understanding of curves in submanifolds, which is
carried out in \S 2. In \S 3 we roughly describe the components of $\s_r(X)$
for an arbitrary variety, and give a precise description for 
$\s_6(\pp 3\times \pp 3\times \pp 3)$, giving a normal
form for a point of each component. In \S 4 we carry out our case by 
case analysis.

\S 5 consists of a treatment of the cases that were overlooked in an
earlier version of this article.

\smallskip

\noindent{\it Acknowledgements.} Supported by NSF grant DMS-0305829.
The key step for this paper was discovered while the   author was a guest 
of L. Manivel at the
University of Grenoble. He  would like to thank
L. Manivel for useful discussions  and providing him a respite from the climate in
Atlanta. He also thanks the anonymous referee who gave a very careful
reading to an earlier version of this paper.

\section{Taylor series for curves on submanifolds}  

Let $V$ be a vector space and $\BP V$ the associated projective
space of lines through the origin in $V$. If $x\in V$ we let
$\hat x\subset V$ denote the line through $x$ and
$[x]\in \BP V$ the corresponding point in projective
space.  If $Z\subset \BP V$
is a set, we let $\hat Z\subset V$ denote the corresponding cone in $V$.

We begin with some very general local differential geometry:

\begin{lemma}\label{curvelemma} Let $X\subset \BP V$  be an analytic submanifold. Let $[x_0]\in X$
and
choose a splitting $V= \hat x_0\oplus T_{[x_0]}X\oplus N_{ [x_0]}X$. (This
is necessary because in affine and projective geometry the tangent and
normal
spaces are only well defined as quotient spaces of $V$.)  
Having done so, we may and will identify abstract and embedded tangent spaces.

Let $x(t)\subset V$ be an analytic curve on $\hat X$ such that $x(0)=x_0$.  
Let $F_j\in S^jT^*_{[x_0]}X\ot N_{[x_0]}X$ denote the $j$-th Fubini
form of $X$ at $[x_0]$ (see \cite{IvL} p. 107 for a definition or one can use the coordinate
definition given in the proof below).
Write, in local coordinates
$$
x(t)= x_0 +tx_1 + t^2 x_2 +\cdots
$$
Then there exists a sequence of elements $y_1,y_2,\hdots \in T_{[x_0]}X$
such that 

i. $x_1=y_1$.  

ii. $x_2= F_2(y_1,y_1) + y_2$.

iii. $x_3=F_3(y_1,y_1,y_1) + F_2(y_1,y_2) + y_3$.

iv. In general, 
$$
x_k=\sum_{j=2}^k \ 
\sum_{\tiny{\begin{matrix}\l_1+2\l_2+\cdots +f\l_f =k \\ \l_1+\cdots +\l_f=j
\end{matrix} }}
 F_j((y_1)^{\l_1},(y_2)^{\l_2}\hd (y_f)^{\l_f}) +y_k .
$$
\end{lemma}

\begin{proof}
First note that despite the choices of splittings,   each term is well defined because
of the lower order terms that appear with it. Also note that ii. is well known in the classical
geometry of surfaces, where $F_2=II$ is the projective {\it second
fundamental form} of $X$ at $[x_0]$.
Take adapted coordinates $( w^{\a},z^{\mu})$ such that $[x_0]=(0,0)$ and
$T_{[x_0]}X$ is spanned by the first $n$ coordinates ($1\leq \a\leq n$). Then locally $X$ is given by equations
\begin{equation}\label{grapheqn}
z^{\mu}=f^{\mu}(w^{\a})
\end{equation} and
$$
F_{k}(\frp {w^{i_1}}\hd \frp{w^{i_k}}) =\sum_{\mu}\frac {\partial^kf^{\mu}}
{\partial w^{i_1} \hd \partial w^{i_k}} \frp {z^{\mu}}
$$
Now write out the vectors $x_j$ in components, substitute into \eqref{grapheqn}
and compare powers of $t$. The result follows.
\end{proof} 

Let $\tau_{k+1}(X)$ denote  the set of all points of the form
$[x_0+x_1+\cdots + x_k]$, where the $x_j$ are as above, so $\tau_1(X)=X$
and $\tau_2(X) $ is the tangential variety of $X$, the union of all
points on all tangent lines to $X$. Note that $\tau_k(X)$ is
not the $k$-th osculating variety for $k>2$, and that the expected dimension
of $\tau_k(X)$ is $kn$. Let $T_{k+1,[x_0]}X$ denote the set of all
points of the form $[x_0+x_1+\cdots + x_k]$ so $\tau_{k+1}(X)=\cup_{[x_0]\in X}
T_{k+1,[x_0]}X$.


\medskip

\section{Components of $\s_k(X)$}
Let $X\subset \BP V$ be a smooth projective variety.
Let $\s_k^0(X)$ denote the set of 
points   in 
$\s_k(X)\backslash \s_{k-1}(X)$ that may be written as the sum of $k$   points on $X$.   
Given $[p]\in \s_k(X)$,  there exist analytic curves $p_1(t)\hd p_k(t)$ in
$\hat X$, such that for $t\neq 0$, $[p(t)]:= [p_1(t)+\cdots +p_k(t)]$ is
in $\s^0_k(X)$ and $[p(0)]=[p]$.  This is true because in the Zariski
topology, a set that is open and dense is also dense in the
classical topology. If one considers the set of honest
secant $\BP^{k-1}$'s as a subset of the Grassmanian of all
$\BP^{k-1}$'s in $\BP V$, this set is Zariski open and dense in
its Zariski closure  
and therefore open and dense as a subset of the
Zariski closure in the analytic topology.

The point $p(0)$ is in the
 limiting $k$-plane corresponding to the first nonvanishing term
in the Taylor series for
$[p_1(t)\ww p_2(t)\ww \cdots \ww p_k(t)]$,  considered as a a point in
the  Grassmannian $G(k,V)\subset \BP (\La k V)$.
  When taking a limit, there will be
  points $q_1\hd q_s$ such that each of the $p_j$'s   limits
to one of the $q_{\a}$'s, e.g., $p_1(t)\hd p_{a_1}(t)$ limits to $q_1$,
$p_{a_1+1}(t)\hd p_{a_2}(t)$ limits to $q_2$ etc... 
We consider separately the limits
$L_1$ of $p_1(t)\ww\cdots \ww p_{a_1}(t)$,
etc...

Now assume that $p_i(0)$  
 is a general point of $X$ for each $i$.
(In particular, this assumption is automatic if $X$ is homogeneous.)
We say the limiting $k$-plane is {\it standard} if
$L_1\ww L_2\ww\cdots \ww L_s\neq 0$, i.e., if we may consider
the limiting linear spaces associated to each of the $q_j$'s separately.
Otherwise we say the limiting $k$-plane is {\it exceptional}.
An example of an exceptional limit is when there are two limiting points $q_1,q_2$, but $q_2$ is in
the tangent space of $q_1$. Another example of an exceptional limit
is if   at least three points go  towards $q_1$,  $q_2 \in T_{3,q_1}X$
and $ T_{3,q_1}X$ is
not the entire ambient space.

\subsection{Standard components for any variety}
We first consider standard limits, so we may restrict our
study to how $a$ curves $p_1(t)\hd p_a(t)$ can limit
to a single point $x$.
Write the Taylor series of $p_j(t)$ as
$$p_j(t)= x + x^j_1t+ x^j_2t^2 + \cdots
$$ 
Consider first for simplicity the case $a=3$. The first
possible nonzero term is

$$
t^2 x\ww (\xxx 21\ww\xxx 3 1 - \xxx 11
\ww \xxx 31 + \xxx 11\ww\xxx 21).
$$
If this term is nonzero then $p\in \s_{k-1}(X)$ because any
point on
such a plane is also on a line of the form $x\ww v$ with
$v\in T_{[x]}X$ and so our three curves only really contribute as two.
An easy exercise shows that if this term is
zero then either all three terms must be
a multiple of one of them, say $\xxx 11$ or two must
be equal, say $\xxx 21=\xxx 31$. But the second case also leads
to $p\in \s_{k-1}(X)$ when we examine the $t^3$ coefficient,
so we must have all a multiple of $\xxx 11$. So the only
type of term we can have is the span of
$$
x, v, II(v,v) + w
$$
where $v,w\in T_{[x]}X$. In other words, under the
hypotheses that $p\notin\s_{k-1}(X)$,
the only possible limit is a point of $\t_3(X)$. 
Similarly for $a=4$, the only possible limit is a point of
$\t_4(X)$.

 A new phenomenon occurs when $a=5$. We  may obtain a point
of $\t_5(X)$ as above, but a second possibility occurs.
Define $\t_5(X)'$ to be the union of all points 
  in the span 
 of 
 $$  x ,x_1,II(x_1,x_1)+x_2,y_1,II(y_1,y_1)+y_2,
 $$
where $x_1,x_2,y_1,y_2\in T_{x}X$.
In the notation of above, this will occur if
$\xxx 11,\xxx 21,\xxx 31,\xxx 41$ span a  plane in $T_{[x]}X$.

For the case where all six points limit to the same point there are
three possibilities. The first case yields
  $\t_6(X)$. For the second, define  $\t_6(X)'$ to be the union of all
points in the span of
$\langle x,x_1,II(x_1,x_1)+x_2,F_3(x_1,x_1,x_1)+II(x_1,x_2)+x_3,y_1,II(y_1,y_1)+y_2\rangle$
where $x_1,x_2,x_3,y_1,y_2\in T_{ x }X$. For the third, define
$\t_6(X)''$ to be the union of all points in the span of 
$\langle x,x_1,II(x_1,x_1)+x_2,F_3(x_1,x_1,x_1)+II(x_1,x_2)+x_3, II(x_1,x_1)+y_2,
F_3(x_1,x_1,x_1)+II(x_1,y_2)+y_3\rangle$, where again, the $x_j$ and
$y_j$ are points of $T_{ x }X$.

\smallskip

We now must take the span of $s$ such points. In general,
given algebraic varieties
$Y_1\hd Y_s\subset\BP V$, define their
{\it join} $J(Y_1\hd Y_s)\subset \BP V$ to be the Zariski
closure of the union of all $\pp{s-1}$'s spanned by points
$y_1\hd y_s$ with $y_j\in Y_j$.

\smallskip

From the above discussion, we obtain all standard components 
of $\s_6(X)\backslash\s_5(X)$ where $X\subset \BP V$ is
  any variety as follows. (Note that these components will not in general be
  disjoint.)

\begin{enumerate}

\item $\s_6^0(X)$,
\item $J(\s_4(x),\t_2(X))$, 
\item $J(\s_3(X),\t_3(X))$,
\item $J(\t_3(X),\t_3(X))$, 
\item $J(\s_2(X),\t_2(X),\t_2(X))$,
\item $J(\t_2(X),\t_2(X),\t_2(X))$,
\item $J(\s_2(X),\t_4(X))$,
\item $J(\t_2(X),\t_4(X))$,
\item $J( X ,\t_2(X),\t_3(X))$
\item $J(X,\t_5(X))$
\item $J(X,\t_5(X)')$  
\item $\t_6(X)$,
\item $\t_6(X)'$
\item $\t_6(X)''$.

\end{enumerate}

\medskip

\subsection{Standard components for $Seg(\BP A\times\BP B\times \BP C)$}
In our
case of  $X=Seg(\BP A\times\BP B\times \BP C)$,
 $X$   is not only homogenous,
but a compact Hermitian symmetric space of rank three. Its only
nonzero differential invariants are
the second fundamental
form $II=F_2$ and the third fundamental form $III$ 
(see \cite{LMsel}, theorem 4.1). The third fundamental
form  is the component of $F_3$ taking image in $N_{x}X/II(S^2T_{[x]}X)$
(see \cite{IvL}, p. 96). 
Unlike the full
  $F_3$, it is a well defined tensor $III\in S^3T^*_xX\ot (N_{x}X/II(S^2T_{[x]}X))$.
Similar to the
situation in  lemma \ref{curvelemma}, we choose a splitting of $N_xX$ to make
$III$ take values in a subspace of $N_{x}X$ instead of a quotient space.

Having such simple differential invariants makes it possible
to have normal forms for elements of each standard component.
We may write any element of
$X$ as $p=[a_1\ot b_1\ot c_1]=[a_1b_1c_1]$,
where all vectors are nonzero. Here and in what follows,
the $a_j$'s are elements of $A$, $b_j$'s of $B$ and $c_j$'s of $C$,
and we omit 
$\otimes$ in the notation   for brevity.  Any element of
$ T_pX$ may be written as
$  
(a_1b_1c_2+  a_1b_2c_1+a_2b_1c_1)$, where
we allow the possibility of some (but not all) of $a_2,b_2,c_2$ to
be zero.  

 If $p=[a_1  b_1  c_1]$
and $v=a_1b_1c_2+  a_1b_2c_1+a_2b_1c_1, w=a_1b_1c_3+  a_1b_3c_1+a_3b_1c_1,
u=a_1b_1c_4+  a_1b_4c_1+a_4b_1c_1
\in T_p X$ then, with the obvious choice of
splitting;
$$
II(v,w)= a_1b_2c_3+a_1b_3c_2+  a_2b_1c_3+  a_3b_1c_2+a_2b_3c_1+a_3b_2c_1 
$$
and
$$
III(u,v,w)= a_4b_2c_3+a_4b_3c_2+  a_2b_4c_3+  a_3b_4c_2+a_2b_3c_4+a_3b_2c_4 .
$$

Here are explicit normal forms for
elements of each standard component when $X=Seg(\BP A\times \BP B\times \BP C)$.
   Repetitions and zeros
among elements of the sets $\{ a_j\}$,
$\{ b_j\}$ and $\{c_j\}$  are   allowed as long as they do not force $p$
into $\s_5(X)$.

\begin{enumerate}
\item  
$p= a_1b_1c_1+\cdots +a_6b_6c_6\in\s_6^0(X)$
\smallskip
\item 
$p= a_1b_1c_1+\cdots +a_4b_4c_4+a_5b_5c_5+  (a_5b_5c_6+a_5b_6c_5+a_6b_5c_5)  \in J(\s_4(x),\t_2(X))$
\smallskip

\item $p= a_1b_1c_1+   a_2b_2c_2+  a_3b_3c_3 + a_4b_4c_4+
 (a_4b_4c_5+a_4b_5c_4+a_5b_4c_4)
+[(a_4b_4c_6+a_4b_6c_4+a_6b_4c_4)+ 2(a_4b_5c_5+a_5b_4c_5+a_5b_5c_4) ]\in J(\s_3(X),\t_3(X)) $
\smallskip

\item $p= a_1b_1c_1+   
(a_1b_1c_2+  a_1b_2c_1+a_2b_1c_1)
+[(a_1b_1c_3+  a_1b_3c_1+a_3b_1c_1)+ 2(a_1b_2c_2+  a_2b_1c_2+a_2b_2c_1)] 
+ a_4b_4c_4+
 (a_4b_4c_5+a_4b_5c_4+a_5b_4c_4)
+[ (a_4b_4c_6+a_4b_6c_4+a_6b_4c_4)+ 2(a_4b_5c_5+a_5b_4c_5+a_5b_5c_4) ]\in J(\t_3(X),\t_3(X)) $

\smallskip

\item $p= a_1b_1c_1+   a_2b_2c_2+  a_3b_3c_3 +  (a_3b_3c_4+a_3b_4c_3+a_4b_3c_3)
+a_5b_5c_5+  (a_5b_5c_6+a_5b_6c_5+a_6b_5c_5)  \in J(\s_2(X),\t_2(X),\t_2(X))$

\smallskip

\item $p= a_1b_1c_1+   (a_1b_1c_2+ a_1b_2c_1+
a_2b_1c_1)+  a_3b_3c_3 +  (a_3b_3c_4+a_3b_4c_3+a_4b_3c_3)
+a_5b_5c_5+  (a_5b_5c_6+a_5b_6c_5+a_6b_5c_5)  \in J(\t_2(X),\t_2(X),\t_2(X))$

\smallskip

\item $p= a_1b_1c_1+   a_2b_2c_2+  a_3b_3c_3 + (a_3b_3c_4+a_3b_4c_3+a_4b_3c_3) 
+[(a_3b_3c_5+a_3b_5c_3+a_5b_3c_3)+ 2(a_3b_4c_4+a_4b_3c_4+a_4b_4c_3) ] 
+ [(a_3b_3c_6+a_3b_6c_3+a_6b_3c_3)+  6a_4b_4c_4 + (a_3b_4c_5+a_3b_5c_4+a_4b_3c_5+a_5b_3c_4+a_4b_5c_3+a_5b_4c_3)]\in J(\s_2(X),\t_4(X))$

\smallskip

\item $p= a_1b_1c_1+ (a_1b_1c_2+  a_1b_2c_1+a_2b_1c_1)+  a_3b_3c_3 + (a_3b_3c_4+a_3b_4c_3+a_4b_3c_3) 
+[(a_3b_3c_5+a_3b_5c_3+a_5b_3c_3)+ 2(a_3b_4c_4+a_4b_3c_4+a_4b_4c_3) ] 
+ [(a_3b_3c_6+a_3b_6c_3+a_6b_3c_3)+  6a_4b_4c_4 + (a_3b_4c_5+a_3b_5c_4+a_4b_3c_5+a_5b_3c_4+a_4b_5c_3+a_5b_4c_3)]\in J(\t_2(X),\t_4(X))$
\smallskip

\item $p= a_1b_1c_1+   a_2b_2c_2+  (a_2b_2c_3+a_2b_3c_2+a_3b_2c_2)    + a_4b_4c_4+
(a_4b_4c_5+a_4b_5c_4+a_5b_4c_4)
+[(a_4b_4c_6+a_4b_6c_4+a_6b_4c_4)+ 2(a_4b_5c_5+a_5b_4c_5+a_5b_5c_4) ]
\in J( X ,\t_2(X),\t_3(X)) $
\smallskip

\item $p= a_1b_1c_1+   
a_2b_2c_2+  (a_2b_2c_3+a_2b_3c_2+a_3b_2c_2)    +  
[(a_2b_2c_4+a_2b_4c_2+a_4b_2c_2) + 2(a_2b_3c_3+a_3b_2c_3+a_3b_3c_2)  
 ]
+
[(a_2b_2c_5+a_2b_5c_2+a_5b_2c_2) +  6a_3b_3c_3+ 
(a_2b_3c_4+a_2b_4c_3+a_3b_2c_4+a_4b_2c_3+a_3b_4c_2+a_4b_3c_2)]
+
[(a_2b_2c_6+a_2b_6c_2+a_6b_2c_2) +  2(a_4b_3c_3+a_3b_4c_3+a_3b_3c_4) + 
(a_2b_3c_5+a_2b_5c_3+a_3b_2c_5+a_5b_2c_3+a_3b_5c_2+a_5b_3c_2)+ 2(a_2b_4c_4+a_4b_2c_4+a_4b_4c_2)]
\in J( X ,\t_5(X) ) $
\smallskip
\item $p= a_1b_1c_1+   
a_2b_2c_2+  (a_2b_2c_3+a_2b_3c_2+a_3b_2c_2)    +  
[(a_2b_2c_4+a_2b_4c_2+a_4b_2c_2) + 2(a_2b_3c_3+a_3b_2c_3+a_3b_3c_2)  
 ]
+
 (a_2b_2c_5+a_2b_5c_2+a_5b_2c_2)    +  
[(a_2b_2c_6+a_2b_6c_2+a_6b_2c_2) + 2(a_2b_5c_5+a_5b_2c_5+a_5b_5c_2)  
 ]
\in J( X ,\t_5(X)' ) $
\smallskip
\item $p= a_1b_1c_1+ (a_1b_1c_2+  a_1b_2c_1+a_2b_1c_1)
+[(a_1b_1c_3+  a_1b_3c_1+a_3b_1c_1)+ 2(a_1b_2c_2+  a_2b_1c_2+a_2b_2c_1)]
+[(a_1b_1c_4+  a_1b_4c_1+a_4b_1c_1)+ 6a_2b_2c_2+
(a_1b_2c_3+a_1b_3c_2+  a_2b_1c_3+  a_3b_1c_2+a_2b_3c_1+a_3b_2c_1)]
+[(a_1b_1c_5+  a_1b_5c_1+a_5b_1c_1)+   2(a_2b_2c_3+a_2b_3c_2+a_3b_2c_2)+
(a_1b_2c_4+a_1b_4c_2+  a_2b_1c_4+  a_4b_1c_2+a_2b_4c_1+a_4b_2c_1)
+ 2(a_1b_3c_3+  a_3b_1c_3+a_3b_3c_1)]
+[
(a_1b_1c_6+  a_1b_6c_1+a_6b_1c_1)
+   
2(a_2b_2c_4+a_2b_4c_2+a_4b_2c_2)+
2(a_2b_3c_3+a_3b_2c_3+a_3b_3c_2)+
(a_1b_2c_5+a_1b_5c_2+  a_2b_1c_5+  a_5b_1c_2+a_2b_5c_1+a_5b_2c_1)
+ (a_1b_3c_4+a_1b_4c_3+  a_3b_1c_4+  a_4b_1c_3+a_3b_4c_1+a_4b_3c_1) ]
\in\t_6(X)$
\smallskip
\item $p= a_1b_1c_1+ 
(a_1b_1c_2+  a_1b_2c_1+a_2b_1c_1)
+[(a_1b_1c_3+  a_1b_3c_1+a_3b_1c_1)+ 2(a_1b_2c_2+  a_2b_1c_2+a_2b_2c_1)]
+[(a_1b_1c_4+  a_1b_4c_1+a_4b_1c_1)+ 6a_2b_2c_2+
(a_1b_2c_3+a_1b_3c_2+  a_2b_1c_3+  a_3b_1c_2+a_2b_3c_1+a_3b_2c_1)]
+ (a_1b_1c_5+  a_1b_5c_1+a_5b_1c_1) 
+[
(a_1b_1c_6+  a_1b_6c_1+a_6b_1c_1)
+2(a_1b_5c_5+  a_5b_1c_5+a_5b_5c_1)]
\in\t_6(X)'$
\smallskip
\item $p= a_1b_1c_1+ 
(a_1b_1c_2+  a_1b_2c_1+a_2b_1c_1)
+[(a_1b_1c_3+  a_1b_3c_1+a_3b_1c_1)+ 2(a_1b_2c_2+  a_2b_1c_2+a_2b_2c_1)]
+[(a_1b_1c_4+  a_1b_4c_1+a_4b_1c_1)+ 6a_2b_2c_2+
(a_1b_2c_3+a_1b_3c_2+  a_2b_1c_3+  a_3b_1c_2+a_2b_3c_1+a_3b_2c_1)]
+ [ (a_1b_1c_5+  a_1b_5c_1+a_5b_1c_1) + 2(a_1b_2c_2+  a_2b_1c_2+a_2b_2c_1)]
+[
(a_1b_1c_6+  a_1b_6c_1+a_6b_1c_1)
+6a_2b_2c_2+
(a_1b_2c_5+a_1b_5c_2+  a_2b_1c_5+  a_5b_1c_2+a_2b_5c_1+a_5b_2c_1)]
\in\t_6(X)''$

\end{enumerate}

\medskip

\subsection{Exceptional components of $\s_6(Seg(\BP A\times \BP B\times \BP C))$}
For an arbitrary variety, writing down all possible
exceptional components is an intractible problem, but, again
by the simplicity of the Segre,
there are only  two exceptional 
components of $\s_6(X)\backslash \s_5(X)$ and both
occur in 
the case of two limit points.

To see this, say there are two exceptional limit points $p,q$
and all six curves limit to   these. {\it A priori} the  possible exceptional
positions of $p$ and $q$ are that $p\in T_qX$ or that $p\in T_{k,q} X$,
for some $k>1$. However, it is easy to see that an element of $T_{k,q}X$ cannot
be decomposable unless either the \lq\lq new\rq\rq\ vector (denoted $y_k$ in the lemma)
is zero, in which case we are reduced to a point of $\s_5(X)$, or
all the Fubini forms occuring in the $k$-th term are zero and $y_k$
is decomposable - but in this case the $k$-th term is just a point
of $T_qX$ and no new phenomenon occurs.

So say  $p\in T_qX$. 
{\it A priori} there could be five different types of limits, depending on the number
of curves that limit to $p$ and the number that limit to $q$. But $p\in T_qX$ implies
$q\in T_pX$, so by symmetry we are   reduced to three cases,
five points limiting to $p$ and one to $q$, four to $p$ and two to $q$, 
and three to each $p$ and $q$. Without loss of generality take
$p=a_1b_1c_1$ and $q=a_1b_1c_2\in T_p X$. If the $x_1$ term in the expansion
for $p$ is not equal to $q$, then nothing new can occur as the
expansions won't interfere with each other (in fact one
ends up with a point of $\s_5(X)$).

In all cases, there is no ambiguity as to which terms in the Taylor expansions for
$p$ and $q$ must contribute, as there
is a unique choice of terms to wedge together that yield a term of lowest order. 

Consider the case of $5$ points limiting to $p$
and one to $q$. In order to get something
new we must have the first tangent vector to $p$ be $a_1b_1c_2$, as then
we can use the first order term in the expansion of $q$.
 In this case we get a point of the following form, where $0$'s have
been included where terms that ordinarily would not be zero are,
e.g., the first zero represents $0=II(a_1b_1c_2,a_1b_1c_2)$.  
Also for clarity a redundant $a_1b_1c_2$ is   included in double
parentheses.

\begin{align*}
x=&a_1b_1c_1 + (a_1b_1c_2) + [ 
(a_1b_1c_3+a_1b_3c_1+ a_3b_1c_1)+0] \\
& +
[ (a_1b_1c_4+a_1b_4c_1+ a_4b_1c_1) +0+ (a_1b_3c_2+ a_3b_1c_2)]\\
&
+[(a_1b_1c_5+a_1b_5c_1+ a_5b_1c_1) +  0+  (a_1b_3c_3+a_3b_1c_3+a_3b_3c_1)
+(a_1b_4c_2+ a_4b_1c_2)]
\\
&
+
 ((a_1b_1c_2))  + (a_1b_1c_6+a_1b_6c_2+a_6b_1c_2).
\end{align*}
Here $x\in \t_5(X)\subset \s_5(X)$, which can
be seen by making the following substitutions:
$\tilde c_5=c_5+c_6$, $\tilde b_4=b_4+b_6$, $\tilde b_5=b_5-b_6$,
 $\tilde a_4=a_4+a_6$, $\tilde a_5=a_5-a_6$.

If the limit to $p$ is in $\t_5(X)'$ we obtain
\begin{align*}
x=&a_1b_1c_1 + (a_1b_1c_2) + [ 
(a_1b_1c_3+a_1b_3c_1+ a_3b_1c_1)+0] \\
& +
  (a_1b_1c_4+a_1b_4c_1+ a_4b_1c_1) \\
&
+[(a_1b_1c_5+a_1b_5c_1+ a_5b_1c_1) +   (a_1b_4c_4+a_4b_1c_4+a_4b_4c_1)]
\\
&
+
 ((a_1b_1c_2))  + (a_1b_1c_6+a_1b_6c_2+a_6b_1c_2).
\end{align*}
Here we may set $\tilde c_5=c_5+c_3$, $\tilde b_5=b_5+b_3$, $\tilde a_5=a_5+a_3$
to see that  $x\in \s_5(X)$.

\smallskip

Still another possibility exists if
in the cases above, the tangent vector to $q$,
namely $(a_1b_1c_6+a_1b_6c_2+a_6b_1c_2)$ already appears
as one of the terms in the expansion for $p$, in which case we get to examine another term in
the Taylor series for $q$. This can happen when the limit to $p$ is 
a point of $\t_5 (X)$
and we have the coincidences $c_6=c_4$, $b_6=b_3$, $a_6=a_3$, $a_4,b_4=0$. Under these circumstances we
get
\begin{align*}
x=&a_1b_1c_1 + (a_1b_1c_2) + [ 
(a_1b_1c_3+a_1b_3c_1+ a_3b_1c_1)+0] \\
& +
[ (a_1b_1c_4 ) +0+ (a_1b_3c_2+ a_3b_1c_2)]\\
&
+[(a_1b_1c_5+a_1b_5c_1+ a_5b_1c_1) +  0+  (a_1b_3c_3+a_3b_1c_3+a_3b_3c_1)
+0]
\\
&
+
 ((a_1b_1c_2))  + ((a_1b_1c_4+a_1b_3c_2+a_3b_1c_2))\\
&
+[(a_7b_1c_2+a_1b_7c_2+a_1b_1c_7)+(a_1b_3c_4+a_3b_1c_4+a_3b_3c_2)].
\end{align*}
Call this case $EX_1$. If the
limit to $p$ is a point of
  $\t_5(X)'$   there is no comparable term to cancel the term in
the expansion for $q$.

One could try to make the $[(a_7b_1c_2+a_1b_7c_2+a_1b_1c_7)+(a_1b_3c_4+a_3b_1c_4+a_3b_3c_2)]$
term also appear in the limit for $p$, but this again forces too much degeneracy.

We have now examined all possibilities for $5$ curves coming together to one point.

Now say the split between $p$ and $q$ is four/two. Again, if $q$ is not
the initial tangent vector to $p$ we can get nothing new, and if
it is, we obtain:
\begin{align*}
x=&a_1b_1c_1 + (a_1b_1c_2) + [ 
(a_1b_1c_3+a_1b_3c_1+ a_3b_1c_1)+0]
\\
& +
[ (a_1b_1c_4+a_1b_4c_1+ a_4b_1c_1)+0 + (a_1b_3c_2+ a_3b_1c_2)]\\
&
+
(( a_1b_1c_2 )) + (a_1b_1c_5+a_1b_5c_2+a_5b_1c_2)
\\
&+
[(a_1b_1c_6+a_1b_6c_2+a_6b_1c_2)+ 2(a_1b_5c_5+a_5b_1c_5+a_5b_5c_2)].
\end{align*}
Call this case $EX_2$.

Now
consider the case where the tangent vector to $q$ occurs in
the expansion for $p$:
\begin{align*}
x=&a_1b_1c_1 + (a_1b_1c_2) + [ 
(a_1b_1c_3+a_1b_3c_1+ a_3b_1c_1)+0]
\\
& +
[ (a_1b_1c_4 )+0 + (a_1b_3c_2+ a_3b_1c_2)]\\
&
+
(( a_1b_1c_2 )) + ((a_1b_1c_3+a_1b_3c_2+a_3b_1c_2))
\\
&+
[(a_1b_1c_5+a_1b_5c_2+a_5b_1c_2)+ 2(a_1b_3c_2+a_3b_1c_4+a_3b_3c_2)]\\
&
+
[a_1b_1c_6+a_1b_6+a_1b_6c_2+a_6b_1c_2+
a_4b_3c_4+a_1b_5c_4+a_5b_1c_4+a_1b_3c_5+a_3b_5c_2+a_5b_3c_2].
\end{align*}
Set $\tilde c_6=c_6+c_4$, $\tilde c_5=c_5+c_2$, and we see that
$x\in \s_5(X)$. Higher order dengerations are similarly eliminated.

Finally consider a $3$-$3$ split. Assuming $q $ is the first tangent vector to $p$ in
the expansion, we get
\begin{align*}
x=&a_1b_1c_1 + ((a_1b_1c_2)) + [ 
(a_1b_1c_3+a_1b_3c_1+ a_3b_1c_1)+0]
\\
&
+[(a_1b_1c_4+a_1b_4c_1+ a_4b_1c_1)+(a_1b_3c_2+a_3b_1c_2)]
\\
&
+
( a_1b_1c_2 ) + (a_1b_1c_5+a_1b_5c_2+a_5b_1c_2)
\\
&+
[(a_1b_1c_6+a_1b_6c_2+a_6b_1c_2)+ 2(a_1b_5c_5+a_5b_1c_5+a_5b_5c_2)]
\end{align*}
Let $\tilde b_3=b_3+b_4$, $\tilde b_6=b_6+\tilde b_3$, 
$\tilde a_3=a_3+a_4$, $\tilde a_6=a_6+\tilde a_3$,
and $\tilde  c_6=c_6+c_4$ and we see $x\in \s_5(X)$.

 Assuming further coincidences
similarly yields nothing new.

\section{Case by case arguments}
We will use variants of 
the proof of Proposition (17.9) in \cite{BCS}, which
is due to Baur. Since there is a misprint
in the proof in \cite{BCS} (a prohibited re-ordering of indices), we reproduce
a proof here.

\begin{theorem}\label{baurthm} If $A$ is a simple $k$-algebra then the rank of the
multiplication operator $M_A$ is $ \geq 2\tdim A-1$.
\end{theorem}

\begin{proof} Let $n:=\tdim A$ and express $M_A$ optimally as
$M_A=\a^1\ot \b^1\ot c_1 +\cdots + \a^r\ot \b^r\ot c_r$ and assume, to obtain
a contradiction, that $r< 2\tdim A-1$. (We switch notation,
working in $A^*\ot B^*\ot C$, using $\a^i$ to
denote elements of $A^*$, and $\b^i$ to denote elements of $B^*=A^*$.)

By reordering if necessary, we may assume $\a^1\hd \a^n$ is a basis of $A^*$.
Let $b\in \langle \b^n, \b^{n+1}\hd\b^r\rangle\upperp$ be a nonzero element
and consider the left ideal $Ab$. We have
$Ab\subseteq \langle c_1\hd c_{n-1}\rangle$. Let $L\supseteq Ab$ be a maximal  left
ideal
containing $Ab$ and let $n-m= \tdim L$.
Since $\tdim Ab\geq m$ (the minimal dimension of an ideal),
 at least $m$ of the $\b^1(b)\hd \b^{n-1}(b)$ are nonzero
(using again the linear independence of $\a^1\hd \a^n$).
Reorder among $1\hd n-1$ such that the first $m$ are nonzero,
  we have 
$\langle c_1\hd  c_m\rangle \subseteq Ab$ (by the linear independence of $\a^1\hd \a^m$).

Note that $\b^m\hd \b^r$ span $A^*$ (otherwise let $z\in \langle \b^m\hd \b^r\rangle\upperp$
and consider the left  ideal $Az$ which is too small). In particular, restricted to $L$,
a subset of them spans $L^*$. We already know $\b^m\mid_L\neq 0$ so we use that as a first basis
vector. 
Let $\b^{i_1}\hd \b^{i_{n-m-1}}$ be a subset 
of $\b^{m+1}\hd \b^r$ such that together with $\b^m$,
when restricted to $L$ they
form a basis of $L^*$. Now let $j_1\hd j_{r-n}$ be a complementary set of
indices such that $ \{ i_1\hd i_{n-m-1}\}\cup 
\{ j_1\hd j_{r-n}\} =\{m+1\hd r\}$ and take $a\in \langle \a^{j_1}\hd \a^{j_{r-n}}\rangle\upperp$
a nonzero element. Consider the
right  ideal $aA\subset \langle c_1\hd c_{m-1},c_m,c_{i_1}\hd c_{i_{n-m-1}}
\rangle$. For any $y\in A$, there exists $w\in L$ such that
$\b^m(y)=\b^m(w)$ and $\b^{i_f}(y)=\b^{i_f}(w)$
for $1\leq f\leq n-m-1$, so $ay-aw\in \langle c_1\hd c_{m-1}\rangle$.
But since $aw$ and the right hand side are both included in $L$, we conclude $aA\subseteq L$ which
is a contradiction as a left ideal cannot contain a nontrivial right ideal.
\end{proof}

We now specialize to matrix multiplication of
$2\times 2$ matrices which we denote by $MMult$. Let
$X=Seg(\pp 3\times \pp 3\times \pp 3)$. The theorem above implies
$MMult\notin \s_6^0(X)$. We now   show it is not in
any of the other possiblilities.  Let $A$ denote the algebra of
two by two matrices.
Note  that 
 all ideals of $A$ must be of dimension two.

\begin{proposition} $Mmult\notin J(\s_4(X),\t_2(X))$,
$J(\s_2(X),\t_2(X),\t_2(X))$ or $J(\t_2(X),\t_2(X),\t_2(X))$.\end{proposition}

\begin{proof}
Say otherwise, that we had an expression respectively
\begin{align*}
 MMult  =& \a^1\b^1c_1+\cdots +\a^5\b^5c_5+
  (\a^5\b^5c_6+ \a^5\b^6c_5 + \a^6\b^5c_5)    
  \\
 MMult  =& \a^1\b^1c_1+  \a^2\b^2c_2   
+  \a^3\b^3c_3+ (\a^3\b^3c_4+ \a^3\b^4c_3 + \a^4\b^3c_3)
+\a^5\b^5c_5\\
&+ (\a^5\b^5c_6+ \a^5\b^6c_5 + \a^6\b^5c_5)  
  \\
 MMult  =&  \a^1\b^1c_1+  (\a^1\b^1c_2+ \a^1\b^2c_1 + \a^2\b^1c_1)   
+  \a^3\b^3c_3+ (\a^3\b^3c_4+ \a^3\b^3c_4 + \a^4\b^3c_3)
\\
&
+\a^5\b^5c_5+ (\a^5\b^5c_6+ \a^5\b^6c_5 + \a^6\b^5c_5)   
\end{align*}
which we refer to as the first, second and third cases.

\smallskip

We first claim that $\a^1,\a^2,\a^3,\a^4$ or $\b^1,\b^2,\b^3,\b^4$
must be linearly independent. Otherwise,  say
both sets were dependent. Let
$a'\in \langle \a^1,\a^2,\a^3,\a^4\rangle\upperp$, $b'\in\langle\b^1,\b^2,\b^3,\b^4\rangle\upperp$,
we have $a'A=Ab'=\langle c_5,c_6\rangle$, a contradiction.
The same conclusion holds for the $4$-ples of vectors
with indices $1,2,5,6$, and those with indices
$3,4,5,6$.
(Here and in all arguments that follow, when we talk about finding
vectors like $a,a',b,b'$ etc..., we mean nonzero vectors.)

\smallskip

We also note various independences among the $c_j$'s.
In case 3, $c_1,c_3,c_5$ must be independent as otherwise
consider $a\in \langle \a^1,\a^3,\a^5\rangle\upperp$,
$b\in \langle \b^1,\b^3,\b^5\rangle\upperp$. We have
$aA\subseteq \langle c_1,c_3,c_5\rangle$ and in fact equality
by the linear dependence, but similarly for $Ab$, thus
$aA=Ab$, a contradiction.
In cases 1 and 2,  $c_3,c_5,c_6$ must be independent as otherwise
in case 1 we could consider $a\in \langle \a^1,\a^2,\a^4\rangle\upperp$
and in case 2,  $a\in \langle \a^1,\a^2,\a^3\rangle\upperp$. In both
cases we get $aA=\langle c_3,c_5,c_6\rangle$ but taking a corresponding
$b$ yields a contradiction as above.

\smallskip

Cases 1 and 2: without loss of generality we assume
$\a^3,\a^4,\a^5,\a^6$ are independent and consider
$b\in \langle \b^1,\b^2,\b^4\rangle\upperp$ in case 1, 
and $b\in \langle \b^1,\b^2,\b^3\rangle\upperp$ in case 2. In both cases
$Ab\subset \langle c_3,c_5,c_6\rangle$.
In case 2 we have the following matrix mapping the coefficients of
$\a^3,\a^5,\a^6$ to the coefficients of $c_3,c_5,c_6$:
$$
Mult(\cdot , b)= 
\begin{pmatrix}
\b^4(b) &0&0\\
0&\b^5(b)+\b^6(b) &\b^5(b)\\
0& \b^5(b) &0
\end{pmatrix}
$$
and in case 1 the matrix is the same except $\b^3(b)$ replaces $\b^4(b)$.
By the linear indepence of
$c_3,c_5,c_6$ and $\a^3,\a^5,\a^6$, the matrix must have rank two.

There are two subcases to consider depending on whether or
not $\b^3(b)=0$ in case 1 (resp. $\b^4(b)=0$ in case 2).

\smallskip

Subcase 1: If $\b^3(b)=0$ (resp. $\b^4(b)=0$)
 then $\b^5(b)\neq 0$ and $Ab=\langle c_5,c_6\rangle$.
We claim that $\b^1,\b^2,\b^4,\b^5$ (resp. $\b^1,\b^2,\b^3,\b^5$) is a basis of $A^*$ as
otherwise let $b'\in \langle \b^1,\b^2,\b^4,\b^5\rangle\upperp$
(resp. $b'\in \langle \b^1,\b^2,\b^3,\b^5\rangle\upperp$).
We would have $Ab'=\langle c_3,c_5\rangle$ which has a one-dimensional
intersection with $Ab$, thus a contradiction.
Thus at least one of $\b^1,\b^2,\b^4$  (resp. $\b^1,\b^2,\b^3$)
 together with $\b^5$ restricted to
$(Ab)^*$ forms a basis of $(Ab)^*$.

Say $\b^1,\b^5$ gives the basis. Then take 
$a\in \langle\a^2,\a^3,\a^4\rangle\upperp$, we have
$aA\subset \langle c_1,c_5,c_6\rangle$. But now for any $y\in A$
there exists $w\in Ab$ such that $\b^1(y)=\b^1(w)$.
Consider 
$ay-aw\subset \langle c_5,c_6\rangle= Ab$. Since
$aw\in Ab$ we conclude $ay\in Ab$ and thus $aA\subseteq Ab$,
a contradiction. The argument is the same if $\b^2,\b^5$ gives
the basis, just change indices. In case 1, the argument is
still the same if $\b^4,\b^5$ gives the basis, again just change
indices.

Say   $\b^3,\b^5$  gives the basis in case 2. Then take $a\in
\langle \a^1,\a^2,\a^3\rangle\upperp$ so $aA\subseteq \langle
c_3,c_5,c_6\rangle$. Now use   $\b^3$  to show
$aA\subseteq \langle c_5,c_6\rangle = Ab$ to again obtain a
contradiction.

\smallskip

Subcase 2: If $\b^3(b)\neq 0$ (resp. $\b^4(b)\neq 0$)
 then $\b^5(b)=0$ and $\b^6(b)\neq 0$ and
$Ab=\langle c_3,c_5\rangle$.
We claim that $\b^1,\b^2,\b^3,\b^4$ is a basis of $A^*$ as
otherwise let $b'\in \langle \b^1,\b^2,\b^3,\b^4\rangle\upperp$
to get a notrivial intersection $Ab\cap Ab'$ and a contradiction.
Thus at least one of $\b^1,\b^2,\b^4$ 
 (resp. $\b^1,\b^2,\b^3$) together with $\b^3$ (resp. $\b^4$)
 restricted to
$(Ab)^*$ forms a basis of $(Ab)^*$. Each case leads to a
contradiction as in the first subcase, finishing the proof for cases
1 and 2.

\smallskip

Case 3: 
 We first claim that $\a^1,\a^3,\a^5$ are linearly independent. Otherwise
 consider $a\in \langle \a^1,\a^3,\a^5,\a^6\rangle\upperp$
 and $a' \in \langle \a^1,\a^3,\a^4,\a^5\rangle\upperp$.
We would have $aA=\langle c_1,c_3\rangle$, $a'A=\langle c_1,c_5\rangle$
with a nontrivial intersection and thus a contradiction.
  
Now consider
$b\in \langle \b^1,\b^3,\b^5\rangle\upperp$, so
$Ab\subset \langle c_1,c_3,c_5\rangle$.
We have the matrix
$$
Mult(\cdot , b)= \begin{pmatrix}
\b^2(b) &0&0\\
0&\b^4(b)  &0\\
0& 0 &\b^6(b)
\end{pmatrix}
$$
which must have rank two by the linear independence of
$c_1,c_3,c_5$ and $\a^1,\a^3,\a^5$.
By symmetry we may assume $\b^6(b)=0$ so $Ab=\langle c_1,c_3\rangle$
and $\b^1,\b^3,\b^5,\b^6$ are linearly dependent.
We claim that $\b^1,\b^2,\b^3,\b^5$ gives a basis of $A^*$ as otherwise
we could find a $b'\in \langle \b^1,\b^2,\b^3,\b^5,\b^6\rangle\upperp$
with $Ab'=\langle c_3\rangle$, a contradiction.
Thus at least one of $\b^1,\b^3,\b^5$ restricted to $(Ab)^*$ together
with $\b^2$ gives a basis of $(Ab)^*$.

Say $\b^1,\b^2$ gives a basis. Then take $a\in \langle
\a^3,\a^5,\a^6\rangle\upperp$ so $aA\subset\langle c_1,c_2,c_3\rangle$.
But $c_2$ appears in $MMult$ with coefficient $\a^1\b^1$ so we may
argue as above to see $aA=\langle c_1,c_3\rangle =Ab$ to obtain
a contradiction.
Similarly for $\b^3,\b^2$,
using $a\in \langle \a^1,\a^5,\a^6\rangle\upperp$
 and the case of $\b^5,\b^2$ using $a\in 
\langle \a^1,\a^3,\a^5\rangle\upperp$. This concludes
 the proof in case 3.

\end{proof} 

\begin{proposition}   $MMult\notin J(\s_3(X),\t_3(X)), J(\t_3(X),\t_3(X)),
J(X,\t_2(X),\t_3(X))$. \end{proposition}
\begin{proof}
Assume otherwise that
\begin{align*}MMult=& \a^1\b^1 c_1+\a^2\b^2 c_2 +\a^3\b^3 c_3
+\a^4\b^4 c_4 +(\a^4\b^4 c_5+ \a^4\b^5 c_4+ \a^5\b^4 c_4)
\\
&
+ (\a^4\b^4 c_6+ \a^4\b^6 c_4+ \a^6\b^4 c_4)+
2(\a^4\b^5 c_5+ \a^5\b^4 c_5+ \a^5\b^5 c_4)
\end{align*}
for the first  case
   and similarly for the other two cases. Let $b\in \langle \b^1,\b^2,\b^3\rangle\upperp$
so $Ab\subseteq \langle c_4,c_5,c_6\rangle$. Note that
$c_4,c_5,c_6$ must be linearly independent, as otherwise take 
$a\in   \langle \a^1,\a^2,\a^3\rangle\upperp$ and $aA=Ab$.

Now consider the linear map $MMult(\cdot, b)$.
Assume for the moment that $\a^4,\a^5,\a^6$ are linearly independent.
With respect to   bases   $c_4,c_5,c_6$ and $\a^4,\a^5,\a^6$, the map
 $MMult(\cdot, b)$ has matrix
$$
\begin{pmatrix}
\b^4(b)+\b^5(b)+\b^6(b) & \b^4(b)+\b^5(b) &\b^4(b) \\
\b^4(b)+\b^5(b) &\b^4(b)&0\\
\b^4(b)&0&0\end{pmatrix}
$$
 It must have two-dimensional image,
but this can occur only if $\b^4(b)=0$, so we conclude $Ab=\langle c_4,c_5\rangle$.
But now take $a\in \langle \a^1,\a^2,\a^3\rangle\upperp$ and the same argument,
assuming $\b^4,\b^5,\b^6$ are linearly independent,
gives $aA=\langle c_4,c_5\rangle$, a contradiction.

Now say $\b^4,\b^5,\b^6$ fail to be linearly independent and consider
the family of ideals $Ab'$ one obtains as $b'$ ranges over $\langle \b^4,\b^5,\b^6\rangle\upperp$.
These ideals must all be contained in $\langle c_1,c_2,c_3\rangle$ hence they must be
constant, otherwise we would have two left ideals with a nontrivial intersection. But that
means either $c_1,c_2,c_3$ fail to be linearly independent, which gives a contradiction
as usual, or $\a^1,\a^2,\a^3$ fail to be linearly independent. But now if
$\a^1,\a^2,\a^3$ fail to be linearly independent, 
consider $a'\in
\langle \a^1,\a^2,\a^3,\a^4\rangle\upperp$, we have $a'A=\langle c_4,c_5\rangle$
equaling $Ab$ above and
giving a contradiction unless $\a^4,\a^5,\a^6$ also fail to be linearly
independent.

Finally, assuming both $\a^4,\a^5,\a^6$ and $\b^4,\b^5,\b^6$ are dependent, 
consider, in the first two cases 
$\tilde a\in \langle \a^3,\a^4,\a^5,\a^6\rangle\upperp$
and $\tilde b\in \langle \b^3,\b^4,\b^5,\b^6\rangle\upperp$ and in the third case
$\tilde a\in \langle \a^2,\a^4,\a^5,\a^6\rangle\upperp$
and $\tilde b\in \langle \b^2,\b^4,\b^5,\b^6\rangle\upperp$. In all cases we have
$\tilde a A=A\tilde b=\langle c_1,c_2\rangle$ and thus a contradiction.
\end{proof}

\begin{proposition}$MMult\notin J(\s_2(X),\t_4(X)), J(\t_2(X),\t_4(X))$.
\end{proposition}

\begin{proof}
We use normal forms as in  \S 3.
Consider $b\in  \langle \b^1,\b^2,\b^3\rangle\upperp$, so $Ab\subset \langle c_3,c_4,c_5\rangle$.
Assuming $\a^3,\a^4,\a^5$ are linearly independent, we get the same
type of  matrix as above and
conclude $\b^4(b)=0$ and $Ab=\langle c_3,c_4\rangle$. If we also assume $\b^3,\b^4,\b^5$
are linearly independent, taking $a\in  \langle \a^1,\a^2,\a^3\rangle\upperp$ we
have $aA=Ab$ a contradiction.

If both $\a^3,\a^4,\a^5$ and $\b^3,\b^4,\b^5$ are linearly dependent, we
may take $a'\in \langle \a^3,\a^4,\a^5,\a^6\rangle\upperp$ and $b'\in 
\langle \b^3,\b^4,\b^5,\b^6\rangle\upperp$ to get $a'A=Ab'=\langle c_1,c_2\rangle$.

So assume $\b^3,\b^4,\b^5$ are linearly dependent and 
$\a^3,\a^4,\a^5$ are independent.
Consider the family of elements $b''\in\langle \b^3,\b^4,\b^5 \rangle\upperp$ so
$Ab''\subset \langle c_1,c_2,c_3\rangle$ this gives a family of left ideals
in $\langle c_1,c_2,c_3\rangle$ but as above there must be just one ideal,
which must be $\langle c_1,c_2\rangle$ as this is what one obtains for the element
in $\langle \b^3,\b^4,\b^5,\b^6 \rangle\upperp$. But this implies $\b^6(b'')=0$
for all $b''\in\langle \b^3,\b^4,\b^5 \rangle\upperp$, i.e.,
$\b^6\subset \langle \b^3,\b^4,\b^5 \rangle$ so now we may take, in the first case
of the proposition
$b_0\in \langle \b^2,\b^3,\b^4,\b^5,\b^6 \rangle\upperp$
   and
$b_0\in \langle \b^1,\b^3,\b^4,\b^5,\b^6 \rangle\upperp$ in the second to obtain a 
one-dimensional left ideal $Ab_0$.
\end{proof}
 
\begin{proposition}$MMult\notin J(X,\t_5(X))$.\end{proposition}
\begin{proof} Let $b\in\langle\b^1,\b^2,\b^3\rangle\upperp$
thus $Ab\subseteq\langle c_2,c_3,c_4\rangle$. As usual we must have
$c_2,c_3,c_4$ linearly independent. Considering the linear map
$MMult(\cdot,b)$, the only $\a^j$'s that arise are $\a^2,\a^3,\a^4$  and, assuming they
are linearly independent, the $3\times 3$ matrix will have zero determinant only if
$\b^4(b)=0$, but then $Ab=\langle c_1,c_2\rangle$.

 Now consider $a\in \langle\a^1,\a^2,\a^3\rangle\upperp$, we see if
$\b^2,\b^3,\b^4$ are linearly independent we have $aA=Ab$ a contradiction.

If say $\b^2,\b^3,\b^4$ fail to be linearly independent, consider
$b'\in \langle \b^2,\b^3,\b^4,\b^5\rangle\upperp$ so $Ab'=\langle c_1,c_2\rangle$.
Now consider
$b''\in \langle \b^2,\b^3,\b^4,\b^6\rangle\upperp$.  We have
$Ab''\subset \langle c_1,c_2,c_3\rangle$ so it has a nonzero intersection with $Ab'$,
hence a contradiction unless it equals $Ab'$. But examining the matrix, even without
assuming independence of $\a^1,\a^2,\a^3$, the only way $Ab''=Ab'$ is if $\b^5\subset
\langle \b^2,\b^3,\b^4,\b^6\rangle$, and if this occurs we simply take
$\tilde b\subset \langle \b^2,\b^3,\b^4,\b^5,\b^6\rangle\upperp$ to obtain a
one dimensional ideal $A\tilde b=\langle c_1\rangle$.
\end{proof}

\begin{proposition} $MMult\notin\t_6(X)$\end{proposition}
\begin{proof}
  Let $b\in\langle \b^1,\b^2,\b^3\rangle\upperp$ and consider
the ideal $Ab\subset\langle c_1,c_2,c_3\rangle$.
Again, $c_1,c_2,c_3$ must
be linearly independent.  If $\a^1,\a^2,\a^3$  
are  also linearly independent, then the only way for $Ab$ to
be two dimensional is if $\b^4(b)=0$.   Thus $Ab=\langle c_1,c_2\rangle$. The same 
reasoning applied to $a\in \langle \a^1,\a^2,\a^3\rangle\upperp$ implies
$aA=\langle c_1,c_2\rangle$, a contradiction, so at least one of the sets
$\a^1,\a^2,\a^3$, $\b^1,\b^2,\b^3$ must be linearly dependent.

Say just one set, e.g., $\b^1,\b^2,\b^3$ fails to be linearly independent,
then by the reasoning above $\b^4\subset\langle \b^1,\b^2,\b^3\rangle$ because
$\b^4(b)=0$ for all $b\in \langle \b^1,\b^2,\b^3\rangle\upperp$. But then
there exists   $b'\in \langle \b^1,\b^2,\b^3,\b^4,\b^5 \rangle\upperp$
yielding an ideal $Ab'=\langle c_2\rangle$, a contradiction.

 Finally say both sets are
dependent.  Consider $b''\in \langle \b^1,\b^2,\b^3,\b^4\rangle\upperp$ and
$a''\in \langle \a^1,\a^2,\a^3,\a^4\rangle\upperp$. We have $Ab''= a'' A=\langle c_1,c_2\rangle$,
a contradiction.
\end{proof}

\begin{proposition} $MMult\notin J(X,\t_5(X)'), \t_6(X)',\t_6(X)'', EX_1,EX_2$\end{proposition}
\begin{proof} 
Consider $b\in\langle \b^1,\b^2,\b^3\rangle\upperp$ in the first four
cases and $b\in\langle \b^1,\b^3,\b^5\rangle\upperp$ for the last.
In all cases 
$Ab$ is a fixed two dimensional ideal (e.g. $\langle c_2,c_5\rangle$ in the first)
but taking  $a \in\langle \a^1,\a^2,\a^3\rangle\upperp$ in the first four
cases and $a\in\langle \a^1,\a^3,\a^5\rangle\upperp$ for the last
yields the same two dimensional ideal, hence a contradiction.
\end{proof}

\section{Cases overlooked in the original article}

Let $X\subset \BP V$ be a projective variety and let $p\in \s_r(X)$. Then
there exist curves $x_1(t)\hd x_r(t)\subset \hat X$ with
$p\in \tlim_{t\ra 0}\langle x_1(t)\hd x_r(t)\rangle $.
We are interested in the case when
$\tdim \langle x_1(0)\hd x_r(0)\rangle <r$. In [1] 
  it was
mistakenly asserted that the only way for this to happen was for
some of the points to coincide. In what follows I show that the cases
I neglected to account for also cannot be matrix multiplication, filling
the gap in the proof.

Use the notation
$x_j=x_j(0)$.    Assume 
for the moment 
that $x_1\hd x_{r-1}$ are linearly independent.
Then we may write $x_r=c_1x_1+\cdots + c_{r-1}x_{r-1}$ for some
constants $c_1\hd c_{r-1}$. Write each curve
$x_j(t)=x_j+tx_j'+t^2x_j'' +\cdots$ where derivatives are taken at $t=0$.

Consider the Taylor series  
\begin{align*}
&x_1(t)\ww\cdots\ww x_r(t) =\\
&(x_1+tx_1'+t^2x_1'' +\cdots)\ww \cdots \ww (x_{r-1}+tx_{r-1}'+t^2x_{r-1}'' +\cdots)
\ww (x_{r}+tx_r'+t^2x_r'' +\cdots) 
\\
&=t((-1)^r(c_1x_1' +\cdots c_{r-1}x_{r-1}'-x_r')\ww x_1\ww\cdots\ww x_{r-1}) + t^2(...) +\cdots
\end{align*}
If the $t$ coefficient is nonzero, then $p$
lies in the linear span of $x_1\hd x_{r-1},(c_1x_1' +\cdots c_{r-1}x_{r-1}'-x_r')$.

If the $t$ coefficient is zero, we have $c_1x_1'+\cdots + c_{r-1}x_{r-1}'-x_r'=
e_1x_1+\cdots e_{r-1}x_{r-1}$ for some constants $e_1\hd e_{r-1}$.
In this case we must examine the $t^2$ coefficient of the expansion. It is
$$ 
(\sum_{k=1}^{r-1} e_kx_k' + \sum_{j=1}^{r-1}c_jx_j''-x_r'')\ww x_1\ww\cdots\ww x_{r-1}
$$
One continues to higher order terms if this is zero.

\smallskip

More generally, 
if only $x_1\hd x_p$ are linearly independent, use index
ranges $1\leq j\leq p$, $p+1\leq s\leq r$ and write
$x_s=c^j_sx_j$ (summation convention), then
  the first possible nonzero term in the Taylor
series is
$$
x_1\wcdots x_p\ww (c^1_{p+1}x_1'+\cdots + c^p_{p+1}x_p'- x_{p+1}')
\wcdots (c^1_{r}x_1'+\cdots + c^p_{r}x_p'- x_{r}')
$$
which (up to a sign) is the coefficient of $t^{r-p}$.

\subsection*{Examples of such limits in the case of the Segre variety}
We look for linear subspaces $L^s\subset A_1\otc A_n$ such
that $\# (\BP L \cap   Seg(\BP A_1\ctimes \BP A_n))>s$.
Let $v_d(\BP V)\subset \BP S^dV$ denote the $d$-th Veronese.
Consider
\begin{equation}\label{space}
\langle v_{d_1}(\pp{i_1-1})\ctimes  v_{d_k}(\pp{i_k-1})\rangle
\end{equation}
where we fix embeddings 
$$S^{d_j}\BC^{i_j}\subset (\BC^{i_j})^{\ot d_j}
\subset A_{j_1}\otc A_{j_{d_j}}
$$

with $d_1+\cdots d_k=n$ and
$$
s\leq\binom{i_1+d_1-1}{d_1}\cdots \binom{i_k+d_k-1}{d_k}
$$
If strict inequality holds, we take a linear section of
\eqref{space} to get an $s$-dimensional linear space.

For example, if we take each $d_j=1$, we are just choosing
linear subspaces $A'_j\subset A_j$ of dimensions $a_j'$
such that $a_1'\cdots a_n'=s$.  

There are other examples, e.g., $\pp 4\cap Seg(\pp 2\times \pp 2)$ is
six ($=deg(Seg(\pp 2\times \pp 2))$) points for a generic $\pp 4$.

\subsection*{Conditions for the  $t$ coefficient  to be
zero in the $3$-factor Segre case}
Write $x_r=\xi^jx_j$, and
assume $x_1\hd x_{r-1}$ are linearly independent. Let $x_{\a}=\tilde a_{\a}\tilde b_{\a}\tilde c_{\a}$ for $1\leq \a\leq r$.
Let $A'=\langle \tilde a_{\a}\rangle$, write $a'=\tdim A'$, and let
$(a_1\hd a_{a'})=(a_i)$ be a basis of $A'$. 
Write $\tilde a_{\a}=\a^i_{\a}a_i$. Do the same for $B',C'$.
We have the   equations
$$
\a^i_r\b^k_r\g^l_r= \xi^{\a}\a^i_{\a}\b^k_{\a}\g^l_{\a}
\ \ \forall i,k,l.
$$

Now assume further   that the $t$ term vanishes, i.e.,
$x_r'=\xi^jx_j'$. Write 
\begin{align*}
x_{\a}'&=a_{\a}'\tilde b_{\a}\tilde c_{\a} +
\tilde a_{\a}  b_{\a}'\tilde c_{\a} + \tilde a_{\a}\tilde b_{\a}c_{\a}'\\
&=\b^k_{\a}\g^l_{\a}a_{\a}'b_kc_l+\a^i_{\a}\g^l_{\a}  a_{i}b_k'c_l
+\a^i_{\a}\b^k_{\a}  a_{i}  b_k'c_l.
\end{align*}
We obtain the  following equations relating $a_j',a_r'$:
\begin{equation}\label{aeqns}
\b^k_r\g^l_ra_r'\equiv  \xi^j\b^i_j\g^l_ja_j' \tmod A'
\end{equation}
which span 
at most an  $N_A:=\tdim \langle \tilde b_1\ot \tilde c_1\hd
\tilde b_r\ot \tilde c_r\rangle$ dimensional space of equations (in particular,
at most $b'c'$).
The $t$-term vanishes iff all these equations and their $B,C$ counterparts
hold. 

If the $t$ term does not vanish we get (up to) $N_A$ new elements
of $A$ appearing beyond the elements of $A'$, and similarly for $B,C$.

If the $t$ term vanishes, then there will be three types of
vectors from $A$ appearing in the $t^2$ term: elements of $A'$,
  the $a_{\a}'$  that live in the solution space to \eqref{aeqns}
and appear as coefficients of $II$,
and the $a_{\a}''$ which may span an $N_A$ dimensional subspace.

\subsection{Some Lemmas}
We assume $a'\leq b'\leq c'$ and use notations as above.

\begin{lemma}\label{outcombos}
Assume $\tdim\langle a_1b_1c_1\hd a_rb_rc_r\rangle <r$.
then

i. If $a'=1$, then $N_A\leq r$.

ii. If $r=4$, then $a'=1$.

iii. If $r=5$, then either $a'=b'=1$ or $a'=b'=c'=2$.

iv. If $r=6$ then either $a'=1$ or $c'\leq 3$.
\end{lemma}
\begin{proof}
All cases are similar, we write out ii. to give a sample calculation.
We assume $a'=b'=c'=2$ and derive a contradiction.
For the first $3$ points there are two possible normal forms:
$a_1b_1c_1, a_2b_2c_2, (\a_1a_1+\a_2a_2)(\b_1b_1+\b_2b_2)(\g_1c_1+\g_2c_2)$
and
$a_1b_1c_1, a_1b_2c_2, a_2(\b_1b_1+\b_2b_2)(\g_1c_1+\g_2c_2)$,
where $a_1,a_2$ are linearly independent elements of $A$, $\a_1,\a_2$
are constants and similarly for $B,C$.
Write an unknown fourth point in the span of the first three
in one of the normal forms as
$(\a_1'a_1+\a_2'a_2)(\b_1'b_1+\b_2'b_2)(\g_1'c_1+\g_2'c_2)$,
expand out in monomials and see that there are no nontrivial
solutions for the fourth point.
\end{proof}
Let
\begin{align*}
&Sub_{a',b',c'}:= 
 \BP \{ T\in A\ot B\ot C\mid \\
& 
\ \ \ \ \ \ \ \ \ \exists A'\subseteq A, B'\subseteq B, C'\subseteq C,
\tdim A'=a',\ \tdim B'=b',\ \tdim C'=c',\ T\in A'\ot B'\ot C'\}
\end{align*}
and let $Sub_{a',b',c'}^0$ be the Zariski open subset of
$Sub_{a',b',c'}$ consisting of elements not in any
$Sub_{a'',b'',c''}$ with $a''\leq a', b''\leq b',c''\leq c'$ and
at least one inequality strict.
Recall that 
$\s_r(Seg(\BP A\times \BP B\times \BP C))\subset
Sub_{r,r,r}$.

\begin{lemma}\label{subspaceprin}
Notations as above.  If $x_1\hd x_r$ are limit points
with $x_r=\xi^jx_j$ such that there
exists curves $x_j(t)$ and
$T\in Sub^0_{r,r,r}\subset A\ot B\ot C$ such that $T$ appears in the
limiting $r$-plane which corresponds to the $t$ coefficient
of the Taylor expansion,  then
for any curves with these limit points such that the $t$ coefficient is zero,
the contribution of $II$ to the $t^2$ coefficient
must be such that $II \equiv 0\tmod A'\ot B'\ot C'$
for some $A',B',C'$, each of dimension $r$ and
$x_j\in A'\ot B'\ot C'$.

More generally, under the same circumstances, if there is
$T$ with
 $T\in Sub_{r-j,rr}^0$,
Then $II\equiv 0\tmod A''\ot B'\ot C'$ where $A'\subset A''$ 
and $\tdim (A''/A')\leq j$.

In particular, if $N_A+a'\geq r$, then $II\equiv 0 \tmod A'\ot B\ot C$
and similarly for permuations.

Similar results hold for higher order invariants.
\end{lemma}

\begin{proof}
We have $r$ curves in each of the vector spaces $A,B,C$ and are examining
their derivatives. The $a''\ot b\ot c + a\ot b''\ot c + a\ot b\ot c''$ type terms
appearing in the $x_j''$ are of the same nature as the terms
appearing in $x_j'$, so they can define a point of $Sub_{r,r,r}^0$.
Thus, were $II\neq 0\tmod A'\ot B'\ot C'$, we could alter our Taylor series to
get a point of $\s_r$ not lying in $Sub_{rrr}$, a contradiction.
 The other assertions are similar.
\end{proof}

\begin{lemma}\label{useful} Let $L^s\subset A\ot B$ be spanned by
$a_1b_1\hd a_sb_s$ and assume all $b_1\hd b_s$ are linearly
independent. Then there is a linear combination of
$a_1b_1\hd a_sb_s$ with all coefficients nonzero that
is in $Seg(\BP A\times \BP B)$ iff $\tdim\langle a_1\hd a_s\rangle =1$.
\end{lemma}

The proof is straightforward. The following lemmas
are consequences of the above lemmas:

\begin{lemma}\label{ttermlem} If $a'=b'=1$
and $r-1=\tdim C'$ then    $II\equiv 0\tmod A'\ot B'\ot C'$.
\end{lemma}

\begin{lemma}\label{a'=1t^2}
If $a'=1$ and there is only  one relation
among   $x_1\hd x_r$, then the terms arising from $II\tmod A'\ot B'\ot C'$ must be 
in $A'\ot B\ot C$.
Moreover, in this situation,  we must have $III\equiv 0\tmod A'\ot B'\ot C' +\tim II$.

If there are $p$ relations among the $x_s$ then at most
$p-1$   vectors which are independent as elements of   $A/ A'$ can appear in $II$
and $III$.
\end{lemma}



\subsection{Matrix multiplication}
Throughout this section we assume
without loss of generality that $a'\leq b'\leq c'$.

  Write the matrix multiplication operator for two by two matrices
as
$M=\phi_1+\phi_2$, and assume $\phi_1\in \s_r$, $\phi_2\in\s_{6-r}$
and $\phi_1$ is not in $\s_r^0$ or any of the cases treated in [1].
We rule out the remaining cases. We must examine the cases $r=3\hd 6$.

\subsection{Case $r=3$.}

Any trisecant  line to the Segre must be in it.
Hence we may assume
$(a',b',c')=(1,1,2)$.  Assume the first nonzero
term is the $t$ coefficient. We have  
$$
p=p_1a_1b_1c_1+p_2a_1b_1c_2+p_3[a_{11}b_1c_1+a_1b_{11}c_1+a_1b_1c_{11}
+a_{12}b_1c_2+a_1b_{12}c_2] 
$$
If any of the triples $a_1,a_{11},a_{12}$, $b_1,b_{11},b_{12}$, $c_1,c_2,c_{11}$ fail
to be linearly independent, we use the $\FS_3$ action to make
that space the target of $M$ (say it is $C$, i.e.,
$M: A^*\times B^*\ra C$).  Take
$b\in\trker\phi_2$ and $a\in \tlker\phi_2$ and consider the
maps $M(\cdot, b)=\phi_1(\cdot,b)$ and $M(a,\cdot)=\phi_1(a,\cdot)$.
Both $M(A^*,b),M(a,B^*)$ are nontrivial ideals, but if $c_1,c_2,c_{11}$
span a two dimensional space the ideals coincide and we have
a contradiction.
So all three of the triples are linearly independent. Now examine the
matrix of the map $M(\cdot, b)$ with respect to the bases
$a_1,a_{11},a_{12}$ and $c_1,c_2,c_{11}$:
$$
\begin{pmatrix}
p_1b_1(b)+p_3b_{11}(b) & p_3 b_1(b) & 0\\
p_2b_1(b)+p_3b_{12}(b) & 0 & p_3b_1(b)\\
p_3b_1(b) & 0 &0\end{pmatrix}
$$
This matrix must have rank two, but that is impossible, a contradiction.
 By lemma \ref{ttermlem}, nothing changes if the $t$ coefficient is zero.

\subsection{Case $r=4$}
  By lemma \ref{outcombos} we must consider   the cases $(1,1,2),(1,1,3)$ and $(1,2,2)$.

\subsubsection{Subcase $(a',b',c')=(1,1,x)$}\label{4,113}
The case $x=2$ is the same as above except take
$b\in  \trker\phi_2\cap  B'\upperp $, $c\in
 \tlker\phi_2\cap C'\upperp$. For the $x=3$ case
   assume for the moment that the $t$ coefficient is nonzero, 
  we have
$$
p\in\langle  a_1b_1c_1,a_1b_1c_2,a_1b_1c_3,
 [a_{11}b_1c_1+a_1b_{11}c_1+a_1b_1c_{11}
+a_{12}b_1c_2+a_1b_{12}c_2+a_{13}b_1c_3+a_1b_{13}c_3 ]
\rangle.
$$
Consider $M:B^*\times C^*\ra A$
  Take
  $b\in \tlker \phi_2\cap b_1\upperp$. Then
$M(b,C^*)=\langle a_1\rangle$ a one dimensional ideal, a contradiction.
Lemma \ref{ttermlem} shows the same argument holds if the $t$ coefficient
is zero.

\subsubsection{Subcase $(a',b',c')=(1,2,2)$}
In this case we have a hyperplane in $\BC^2\ot \BC^2$, thus it intersects
the smooth quadric hypersurface $\pp 1\times \pp 1$ in a plane conic. Up to
equivalence, there are two cases, depending on whether the conic is
smooth or the union of two lines. The second case reduces to (a sum of
cases with) $b'=1$.  
Assuming that no three points are colinear,  we have the normal
form 
$x_1=a_1b_1c_1$, $x_2=a_1b_2c_2$, $x_3=a_1(b_1+b_2)(c_1+c_2)$,
$x_4=a_1(\b^1b_1+\b^2b_2)(\g^1c_1+\g^2c_2)$ with 
$\b^1\g^2=\b^2\g^1$. 
The $t$ term is $x_1\ww x_2\ww x_3$ wedged against
\begin{align*}
&a_4'(\b^1b_1+\b^2b_2)(\g^1c_1+\g^2c_2)+a_1b_4'(\g^1c_1+\g^2c_2)
+a_1(\b^1b_1+\b^2b_2)c_4'\\
&-\b^1\g^2[ a_3'(b_1+b_2)(c_1+c_2)+a_1b_3'(c_1+c_2)+a_1(b_1+b_2)c_3']\\
&-\b^1(\g^1-\g^2)[a_1'b_1c_1+a_1b_1'c_1+a_1b_2c_1']
-\b^2(\g^2-\g^1)[a_2'b_2c_2+a_2b_2'c_2+a_2b_2c_2']\\
&=
\tilde a_1'b_1c_1+\tilde a_2'b_2c_2+\tilde a_3'(b_1c_2+b_2c_1)
+
a_1\tilde b_1'c_1+a_1\tilde b_2'c_2
+
a_1b_1\tilde c_1' + a_1b_2\tilde c_2'
\end{align*}
where $\tilde a_1'=\b^1\g^1a_4'-\b^1\g^2a_3'-\b^1(\g^1-\g^2)a_1'$
etc...

If we take  a nonzero $b\in\tlker \phi_2$, the map $M(\cdot, b): A^*\ra C$
must have rank four or two. We may require additionally that
$\tilde b_2'(b)=0$, in which case $M(\cdot, b)$ has rank at most
three and therefore two. But then  the only way to
obtain a map of rank two is for  $\tilde b_1'(b)=0$, and we
have a nontrivial left ideal $M(A^*,b)=C'$.
But now consider $a\in\trker (\phi_2,a_1)$. $M(a,B^*)=C'$
and we have a contradiction.

From this normal form one can check directly that
$II\equiv 0 \tmod A'\ot B'\ot C'$ so there is no need to
consider the case where the $t$ coefficient is zero.
 
\subsection{Case $r=5$} By lemma \ref{outcombos} we are reduced to
$(1,1,x),(1,2,2),(2,2,2),(1,2,3),(1,3,3)$.

\smallskip

\begin{remark}\label{remcase}
If $a',b',c'\leq 2$, then
$M$ cannot correspond to using the $t$ term in $\phi_1$.
To see  this, we may take $a\in A^*,b\in B^*$ such that
$M(a,B^*),M(A^*,b)\subseteq C'$ by   having $a$ annhilate
$\phi_2$ and $A'$ and similarly for $b$.
Since $\tdim C'=2$ we must have $M(a,B^*)=M(A^*,b)=C'$, 
  a contradiction.
\end{remark}

\subsubsection{Subcase $(a',b',c')=(1,1,x)$}
Argue as in    \S\ref{4,113}. (At most one new vector from each 
of $A,B$ can appear at any given level over the $r=4$ case but
we are allowed to annhilate one new vector here as $\phi_2$ is 
smaller.)

\subsubsection{Subcase
$(a',b',c')=(1,2,2)$}\label{5,122} Consider the $t^2$ term.
Let $\tilde a'$ be the
unique vector mod $A'$ that appears from one derivative. 
Consider $M: B^*\times C^*\ra A$.
Taking $b\in   B'\upperp\cap \tlker\phi_2 $,
we have $M(b,C^*)\subseteq \langle A',\tilde a'\rangle$, but
taking $c\in   C'\upperp \cap \trker\phi_2$,
we have $M(B^*,c)\subseteq \langle A',\tilde a'\rangle$,
a contradiction. The $t^3$ case follows as $III\equiv 0\tmod
\langle A',\tilde a'\rangle \ot B\ot C$ so we may use the same argument.

\subsubsection{Subcase $(a',b',c')=(2,2,2)$}
$N_A,N_B,N_C$ each is $3$ or $4$  and in any of these cases
lemma \ref{subspaceprin} shows $II\equiv 0\tmod A'\ot B'\ot C'$.

\subsubsection{Subcases $(a',b',c')=(1,2,3),(1,3,3)$}\label{5,133}
Let $\phi_2=a_6b_6c_6$. Then taking $b\in B'\upperp$,
$c\in C'\upperp$ we have $M(b,C^*)=M(B^*,c)=\langle A', a_6\rangle$
by lemma \ref{a'=1t^2}.

\subsection{Case $r=6$}
The same argument as in remark \ref{remcase} shows
that if $a',b',c'\leq 3$ with at least one of them at most two,
then $M$ cannot be obtained by using the $t$ term in $\phi_1$.

\subsubsection{Subcases $(a',b',c')=(1,2,2), (1,2,3),(1,3,3)$}  Same proof as 
\ref{5,122} works, for the last two we are allowed   additional vectors
in $B',C'$ as there is no $\phi_2$.

\subsubsection{Subcase: $(a',b',c')=(1,1,x)$}
Argue as in \S\ref{4,113}.

\subsubsection{Subcase: $(a',b',c')=(2,2,2)$}
Here if $II$   is nonzero modulo
$A'\ot B'\ot C'$ there is at most one new vector
in each space appearing at the level of first derivatives
by lemma \ref{subspaceprin}. Call these vectors $\tilde a',
\tilde b',\tilde c'$. Then take
$b\in \langle B',\tilde b'\rangle\upperp$,
$c\in \langle C',\tilde c'\rangle\upperp$, we
have $M(b,C^*)=M(B^*,c)=A'$.
$III$ cannot be nonzero modulo $A'\ot B'\ot C'$ by
lemma \ref{subspaceprin}.

\bibliographystyle{amsplain}
\bibliography{Lmatrix}

\medskip

Address:

Dept. of Mathematics

Texas A\& M University

College Station, TX 77843-3368

\end{document}